\documentclass[12pt]{amsart}
 \newtheorem{thm}{Theorem}[section]
 \newtheorem{cor}[thm]{Corollary}
 \newtheorem{lem}[thm]{Lemma}
 \newtheorem{prop}[thm]{Proposition}
 \theoremstyle{definition}
 \newtheorem{defn}[thm]{Definition}
 \theoremstyle{remark}
 
 \numberwithin{equation}{section}

 \newcommand{\norm}[1]{\left\Vert#1\right\Vert}
 
 \newcommand{\C}{\mathbb{C}}

\begin{document}

\title[A Shift Operator on $L(H^2)$]
 {A Shift Operator on $L(H^2)$}

\author{ Yun-Su Kim }

\address{Department of Mathematics, Indiana University, Bloomington,
Indiana, U.S.A. }

\email{kimys@indiana.edu}


\dedicatory{}

\commby{Daniel J. Rudolph}


\begin{abstract}We give definitions and some properties of the shift operator $S_{L(H^2)}$ and
multiplication operator on $L(H^2)$. In addition, we obtain some
properties of the commutant of the shift operator $S_{L(H^2)}$ and
characterize $S_{L(H^2)}$-invariant subspaces.
\end{abstract}

\maketitle

\section*{Introduction}
For the Hardy space $H^2$ of the open unit disk $\textbf{D}$ and a
holomorphic bounded function $\phi\in{H^{\infty}}$, the shift
operator $S:H^{2}\rightarrow{H^{2}}$ defined by $(Sf)(z)=zf(z)$
for $f\in{H^{2}}$ and the Toeplitz operator
$T_{\phi}:H^{2}\rightarrow{H^{2}}$ defined by
$(T_{\phi}f)(z)=\phi(z)f(z)$ for $z\in\textbf{D}$ and
$f\in{H^{2}}$ are very important operators on the Hardy spaces,
and a classical result, proved by Beurling \cite{1}, links inner
functions with the theory of invariant subspaces for $S$. The
Beurling's Theorem says that every $S$-invariant subspace of $H^2$
is generated by an inner function. In this paper, we provide a
sufficient condition of an invariant subspace for the Toeplitz
operator $T_{\phi}\in{L(H^2)}$ in order to be generated by an
inner function [Theorem \ref{51}].

In section 2, we introduce a new shift operator
$S_{L(H^2)}:L(H^2)\rightarrow {L(H^2)}$ as follows :
\begin{center}$S_{L(H^2)}(T)=ST\texttt{ for }T\in{L(H^2)},$\end{center}
and a multiplication operator $M_{\theta}:L(H^2)\rightarrow
{L(H^2)}$ for $\theta\in{H^\infty}$ defined by
\begin{center}$M_{\theta}T=T_{\theta}{T}$ where $T\in{L(H^2)}$.\end{center} By using
fundamental properties of $S\in{L(H^2)}$ and
$T_{\theta}\in{L(H^2)}$, we also prove same results for the
operators $S_{L(H^2)}$ and $M_{\theta}$ in Corollary \ref{57} and
Lemma \ref{58}. To find some characterizations of the commutant of
$S_{L(H^2)}$, we use the well known fact that the commutant of $S$
is the set $\{T_{\phi}:\phi\in{H^\infty}\}$ in Theorem \ref{71}
and Corollary \ref{83}.

Whenever a new shift operator was introduced, mathematicians gave
an analogue of Beurling's theorem, as examples M. B. Abrahamse and
R.G. Douglas \cite{2}, Hasumi \cite{12}, Sarason \cite{17}, and
Voivhick \cite{19}. In this paper, we also characterize
$S_{L(H^{2})}$-invariant subspaces in Theorem \ref{30} as
following : \vskip0.2cm
 $M$ is a closed invariant subspace for $S_{L(H^{2})}$ if and only
 if there is an inner function $\theta\in{H^{\infty}}$ such that
\begin{center}$M=T_{\theta}M^{\prime}=\{T_{\theta}T^{\prime}:T^{\prime}\in{M^{\prime}}\},$\end{center} where $M^{\prime}$ is a closed
$S_{L(H^{2})}$-invariant subspace of $L(H^{2})$ and $H^{2}$ is
generated by
$\bigcup_{T\in{M^{\prime}}}\{f\in{H^{2}}:f\in{\texttt{range of
}T}\}$.

\section{Preliminaries and Notation}
In this paper, $\C$, ${M}^{-}$, and $L(H)$ denote the set of
complex numbers, the (norm) closure of a set $M$, and the set of
all bounded linear operators from $H$ to $H$ where $H$ is a
Hilbert space, respectively.

\subsection{Hardy Spaces.} \label{1}We refer to \cite{R} for
the basic facts about the Hardy spaces, and recall here the basic
definitions. Let $\textbf{D}$ be the open unit disc.
\begin{defn}\label{60}

The space ${\: H^{2}}$ is defined to be the space of analytic
functions $f$ on $\textbf{D}$ such that the subharmonic function
$|f|^{2}$ has a harmonic majorant on $\textbf{D}$. For a fixed
$z_{0}$ $\in{\textbf{D}}$, there is a norm on $H^{2}$ defined by

      $\norm{f}_{2}=$ inf $\{ u(z_{0})^{1/2}$: $u$ is a harmonic majorant of
      $|f|^{2}\}$.
\end{defn}

We know that $H^{2}$ is also the set of all analytic functions
$f:\textbf{D}\rightarrow{{\C}}$ with power series expansion
$f(z)=\sum_{n=0}^{\infty}a_{n}z^{n}$ for $z\in{\textbf{D}}$ such
that $\sum_{n=0}^{\infty}$$\mid$${a_{n}}$$\mid^{2}$$<\infty$. The
space $H^{2}$ is a Hilbert space with the inner product
$(f,g)=\sum_{n=0}^{\infty}(a_{n},b_{n})=\sum_{n=0}^{\infty}a_{n}\overline{b_{n}}$
for $f(z)=\sum_{n=0}^{\infty}a_{n}z^{n}\in{H^{2}}$ and
$g(z)=\sum_{n=0}^{\infty}b_{n}z^{n}\in{H^{2}}$. Let $m$ be a
harmonic measure for the point $z_{0}$, let
 $L^{2}(\partial{\textbf{D}})$ be the $L^{2}$-space of complex valued
 functions on the boundary of $\textbf{D}$ defined with respect to $m$,
 and let $H^{2}(\partial{\textbf{D}})$ be the set of
 functions $f$ in $L^{2}(\partial{\textbf{D}})$ such that
 $\int_{\partial{\textbf{D}}} f(z)g(z) dz$ = 0 for every $g$ that is
 analytic in a neighborhood of the closure of $\textbf{D}$.
 If $f$ is in $H^{2}$, then there is a function $f^{\ast}$ in
 $H^{2}(\partial{\textbf{D}})$ such that $f({z})$ approaches
 $f^{\ast}(\lambda_{0})$ as $z$ approaches $\lambda_{0}$
 nontangentially, for almost every $\lambda_{0}$ relative to $m$. The map
 $f\rightarrow{f^{\ast}}$ is an isometry from $H^{2}$ onto
 $H^{2}(\partial{\textbf{D}})$. In this way, $H^{2}$ can be
 viewed as a closed subspace of $L^{2}(\partial{\textbf{D}})$.

 We denote by $H^{\infty}$ the Banach space of all bounded analytic
functions $\phi:\textbf{D}\rightarrow{\C}$ with the norm
$\norm\phi_{\infty}=$ sup
$\{\mid$${\phi(z)}$$\mid:z\in{\textbf{D}}\}$. $H^{\infty}$ is a
closed
 subspace in the weak$^{*}$ topology of $L^{\infty}$, in fact, the mapping
 $\phi\rightarrow{\phi^{\ast}}$ is an isometry of $H^{\infty}$
 onto a weak$^{*}$-closed subalgebra of
 $L^{\infty}(\partial\textbf{D})$.
For every $\phi\in{H^{\infty}}$, we can define the Toeplitz
operator $T_{\phi}\in{L(H^{2})}$ by $(T_{\phi}f)(z)=\phi(z)f(z)$
for $z\in\textbf{D}$ and $f\in{H^{2}}$.
\begin{defn}\label{50}
A function $\phi\in{H^{\infty}}$ is said to be \emph{inner} if
$\mid$$\phi(e^{it})$$\mid$$=1$ almost everywhere on
$\partial\textbf{D}$.
\end{defn}
\subsection{Invariant subspaces of Toeplitz operators.}\label{55}
One of the useful operators on a Hilbert space is the shift
operator. First, the shift operator $S\in{L(H^{2})}$ is defined by
$(Sf)(z)=zf(z)$ for $f\in{H^{2}}$. For this shift operator $S$, we
have a useful theorem. A. Beurling's Theorem states every
nontrivial closed $S$-invariant subspace, say \emph{N}, of $H^{2}$
is generated by an inner function, i.e. there is an inner function
$\theta\in{H^{\infty}}$ such that $\emph{ N}=\theta{H^{2}}$. Now
we can address the $T_\phi$-invariant subspace of $H^{2}$ which is
generated by an inner function.
\begin{thm}\label{51}
Let $\phi\in{H^{\infty}}$ such that $\norm{\phi}_{\infty}\leq{1}$.
Let $Y$ be a nontrivial closed $T_{\phi}$-invariant subspace of
$H^{2}$ such that \begin{equation}\label{32}f\in{Y} \texttt{ if
and only if }f\circ{\phi}\in{Y}.\end{equation} Then there is an
inner function $\theta\in{H^{\infty}}$ such that
$Y=\theta{H^{2}}$.
\end{thm}
\begin{proof}
By assumption (\ref{32}), for any $f\in{Y}$, $f\circ{\phi}\in{Y}$.
Since $Y$ is $T_{\phi}$-invariant,
\begin{equation}\label{80}T_{\phi}(f\circ\phi)=\phi\cdot{f(\phi)}\in{Y}\end{equation} where
$\phi\cdot{f(\phi)}$ means that
$[\phi\cdot{f(\phi)}](z)=\phi(z)f(\phi(z))$ for
$z\in{\textbf{D}}$. If $h=Sf$ , i.e., $h(z)=zf(z)$, then
$h\circ\phi=\phi\cdot{f(\phi)}$ and so by equation \eqref{80},
$h\circ\phi\in{Y}$. It follows from assumption (\ref{32}) that
$h=Sf\in{Y}$. Since $f\in{Y}$ is arbitrary, we can conclude that
$Y$ ia a closed $S$-invariant subspace of $H^{2}$. By the
Beurling's Theorem, $Y=\theta{H^{2}}$ for some inner function
$\theta\in{H^{\infty}}$.
\end{proof}
Now let's think about an example of Theorem \ref{51}.
\begin{lem}\label{52}
Let $\phi\in{H^{\infty}}$ be a function defined by $\phi(z)=z^2$.
If $Y=\{h\in{H^{2}}:h(z)=zf(z)$ for some $f\in{H^{2}}\}$, then
$g\in{Y}$ if and only if $g\circ\phi\in{Y}$.
\end{lem}
\begin{proof}
Let $\theta\in{H^{\infty}}$ be an inner function defined by
$\theta(z)=z$ for any $z\in{\textbf{D}}$. Since $Y=\theta{H}^{2}$,
$Y$ is closed. If $g\in{Y}$, i.e. $g(z)=zf(z)$ for some
$f\in{H^{2}}$, then clearly $g\circ\phi\in{Y}$. Conversely if
$g(z)=\sum_{n=0}^{\infty}a_{n}z^{n}\in{H^{2}}$ such that
$g\circ\phi\in{Y}$, then
$(g\circ\phi)(z)=\sum_{n=0}^{\infty}a_{n}z^{2n}\in{Y}$. By the
definition of $Y$, $a_0=0$. Thus $g\in{Y}$.
\end{proof}

\section{Commutant of the Shift Operator $S_{L(H^{2})}$}
\subsection{A Shift Operator on $L(H^{2})$}\label{2}
We know several shift operators defined on some Hilbert spaces,
for examples the shift operators on $H^{2}$, $\emph{l}^{2}$ and so
on. We can try to think about a shift operator defined on a Banach
space which is not a Hilbert space. In this paper, we will define
a shift operator $S_{L(H^2)}:L(H^2)\rightarrow{L(H^2)}$. Note that
$L(H^2)$ is a Banach space which is not a Hilbert space.
\begin{defn}\label{54}
We define the \emph{shift operator}
$S_{L(H^2)}:L(H^2)\rightarrow{L(H^2)}$ by $S_{L(H^2)}(T)=ST$ for
$T\in{L(H^2)}$ where $S\in{L(H^2)}$ is the shift operator defined
in the Section \ref{55}.
\end{defn}
We will prove that $S_{L(H^2)}\in{L(L(H^2))}$ [See Theorem
\ref{56}]. By the definition of $S\in{L(H^2)}$,
$[(S_{L(H^2)}T)f](z)=[S(T{f})](z)=z(T{f})(z)$ for $f\in{H^{2}}$
and $T\in{L(H^2)}$.
\begin{thm}\label{56}
 $S_{L(H^2)}\in{L(L(H^2))}$ and $\norm{S_{L(H^2)}}=1$.
\end{thm}
\begin{proof}
For $a_{i}\in\C(i=1,2)$ and $T_{i}\in{L(H^2)}(i=1,2)$,
$S_{L(H^2)}(a_{1}T_{1}+a_{2}T_{2})=S(a_{1}T_{1}+a_{2}T_{2})
=a_{1}(ST_1)+a_{2}(ST_2)
=a_{1}S_{L(H^2)}(T_1)+a_{2}S_{L(H^2)}(T_2)$ is clear.

For $f\in{H^{2}}$ and $T\in{L(H^2)}$, since $S$ is an isometry on
$H^2$,
$\norm{(S_{L(H^2)}T)f}_{2}=\norm{S(T{f})}_{2}=\norm{T{f}}_{2}
\leq\norm{T}\norm{f}_{2}$. Thus $\norm{S_{L(H^2)}T}\leq\norm{T}$
for any $T\in{L(H^2)}$ and so $\norm{S_{L(H^2)}}\leq{1}$. For
$T=I_{H^2}$(which is the identity operator on $H^{2}$),
$S_{L(H^2)}T=S$, and so $\norm{S_{L(H^2)}T}=\norm{S}=1$. It
follows that $\norm{S_{L(H^2)}}=1$.
\end{proof}
When we use the fact that $S\in{L(H^2)}$ is an isometry, we get
the following result. This result could be an another proof of
Theorem \ref{56}.
\begin{cor}\label{57}
$S_{L(H^2)}$ is an isometry.
\end{cor}
\begin{proof}
For $T\in{L(H^2)}$, $\norm{T}=$
sup$_{f\in{H^{2}}}\frac{\norm{T{f}}_{2}}{\norm{f}_{2}}=$
sup$_{f\in{H^{2}}}\frac{\norm{S(T{f})}_{2}}{\norm{f}_{2}}=$
sup$_{f\in{H^{2}}}\frac{\norm{(S_{L(H^2)}T){f}}_{2}}{\norm{f}_{2}}=\norm{S_{L(H^2)}T}$.
\end{proof}
\subsection{Multiplication Operator on $L(H^2)$}\label{77}
If $(X,\varsigma,\mu)$ is a probability space, for
$\xi\in{L^{\infty}(\mu)}$ we know a multiplication operator
$M_{\xi}$ defined on $L^{2}(\mu)$ such that $M_{\xi}f=\xi{f}$ for
$f\in{L^{2}(\mu)}$ where $\xi{f}$ denotes the pointwise product,
and the operator $T_{\phi}$ defined in Section \ref{1} is also a
multiplication operator defined on $H^{2}$. In the definition of
the multiplication operator on $L(H^2)$, we will use the operator
$T_{\phi}$ as following:
\begin{defn}For $\theta\in{H^\infty}$,
define a \emph{multiplication operator}
$M_{\theta}:L(H^2)\rightarrow{L(H^2)}$ by
$M_{\theta}T=T_{\theta}{T}$ where $T\in{L(H^2)}$.
\end{defn}
Then for $\theta\in{H^\infty}$ and $T\in{L(H^2)}$,
$[(M_{\theta}T)f](z)=\theta(z)(T{f})(z)$ where $f\in{H^{2}}$. You
can check easily $\norm{T_{\theta}}=\norm{\theta}_{\infty}$
\cite{B1}.
\begin{lem}\label{58}
For $\theta\in{H^\infty}$, $M_{\theta}\in{L(L(H^2))}$ and
$\norm{M_{\theta}}=\norm{\theta}_{\infty}$.
\end{lem}
\begin{proof}
For $a_{i}\in\C(i=1,2)$ and $T_{i}\in{L(H^2)}(i=1,2)$,
$M_{\theta}(a_{1}T_{1}+a_{2}T_{2})=a_{1}M_{\theta}T_{1}+a_{2}M_{\theta}T_{2}$
is clear. Let $T\in{L(H^2)}$. Since
$\norm{M_{\theta}T}=\norm{T_{\theta}{T}}\leq\norm{T_{\theta}}\norm{T}
=\norm{\theta}_{\infty}\norm{T}$,
$\norm{M_{\theta}}\leq\norm{\theta}_{\infty}$. Let $T=I_{H^2}$.
Then $M_{\theta}T=T_{\theta}{I_{H^2}}=T_{\theta}$ and so
$\norm{M_{\theta}I_{H^2}}=\norm{T_\theta}=\norm{\theta}_{\infty}$.
Thus $\norm{M_{\theta}}=\norm{\theta}_{\infty}$.
\end{proof}
\subsection{Commutant of the Shift Operators}\label{21}
For $a\in{\C}$, let $u_{a}\in{H^2}$ be the constant function
defined by $u_{a}(z)=a$ for $z\in{\textbf{D}}$. If $f\in{H^2}$
with a power series expansion $f(z)=\sum_{n=0}^{\infty}a_{n}z^{n}$
for $z\in{\textbf{D}}$, then
$f\equiv\sum_{n=0}^{\infty}S^{n}u_{a_n}$, because
$(\sum_{n=0}^{\infty}S^{n}u_{a_n})(z)=\sum_{n=0}^{\infty}z^{n}(u_{a_n}(z))
=\sum_{n=0}^{\infty}a_{n}z^{n}$. We know that the following
Proposition is true (\cite{H}).

\begin{prop}\label{70}
If $T\in{L(H^2)}$ and $ST=TS$, then there is a function
$\theta\in{H^\infty}$ such that $(Tg)(z)=\theta(z)g(z)$ for
$g\in{H^2}$, that is, $T=T_\theta$.
\end{prop}

Recall that the \emph{commutant} of an operator $T\in{L(H)}$ is
the set $\{T_{1}\in{L(H)}:TT_{1}=T_{1}T\}$ and is denoted by
$\{T\}^{\prime}$. By using Proposition \ref{70}, we provide the
maximal abelian subalgebra of $L(H^{2})$ in the following
Corollary.
\begin{cor}\label{81}
The algebra $A=\{T_{\theta}\in{L(H^{2})}:\theta\in{H^\infty}\}$ is
maximal abelian.
\end{cor}
\begin{proof}
Suppose
$T\in{A^{\prime}}=\{T_{1}\in{L(H^{2})}:T_{2}T_{1}=T_{1}T_{2}\texttt{
for any }T_{2}\in{A}\}$. Then clearly $TS=ST$, since $S\in{A}$. By
Proposition \ref{70}, there is a function $\theta\in{H^\infty}$
such that $T=T_\theta$. Clearly $A$ is commutative. It follows
that $A$ is maximal abelian.
\end{proof}
Now we will explore the commutant of the shift operators on
$L(H^2)$. First, if $X\equiv{M_{\theta}}$ for
$\theta\in{H^\infty}$, then
$[((S_{L(H^2)}X)T)f](z)=[((XS_{L(H^2)})T)f](z)=z\theta(z)(Tf)(z)$
for $T\in{L(H^2)}$, $f\in{H^2}$ and $z\in{\textbf{D}}$. Thus
$M_{\theta}\in{\{S_{L(H^2)}\}^{\prime}}$.
\begin{thm}\label{71}
Let $N\subset{L(H^2)}$ be the (norm) closed subspace generated by
$\{S^n:n=0,1,2,\cdot\cdot\cdot\}$ where $S^0=I_{H^2}$ which is the
identity function on $H^2$. Then for every operator
$X\in{\{S_{L(H^2)}\}^{\prime}}$ such that
$X(I_{H^2})\in{\{S\}^\prime}$, there is a function
$\theta\in{H^\infty}$ such that $X\equiv{M_{\theta}}$ on $N$.
\end{thm}
\begin{proof}
Suppose $X\in{\{S_{L(H^2)}\}^{\prime}}$. Since
$X(I_{H^2})\in{\{S\}^\prime}$, by Proposition \ref{70}, there is a
function $\theta\in{H^\infty}$ such that $X(I_{H^2})=T_\theta$.
Since $(S_{L(H^2)}{X})T=(X{S_{L(H^2)}})T$ for any $T\in{L(H^2)}$,
\begin{equation}\label{72}S(X(T))=X(ST).\end{equation}
Let $T=I_{H^2}$. Then by equation \eqref{72},
\begin{equation}\label{78}S{X(I_{H^2})}=X(SI_{H^2})=X(S),\end{equation} and if $T=S$, then
from \eqref{72} and \eqref{78} we have
\[X(S^2)=SX(S)=S(S{X(I_{H^2})}) =S^{2}X(I_{H^2}).\] By the same
way as above, we have
\begin{equation}\label{73}X(S^n)=S^{n}X(I_{H^2}).\end{equation}
for any $n=0,1,2,\cdot\cdot\cdot$.

Since $X(I_{H^2})=T_\theta$, by equation \eqref{73}, $
X(S^n)=S^{n}T_{\theta}=T_{\theta}S^{n}=M_{\theta}(S^{n})$ for any
$n=0,1,2,\cdot\cdot\cdot$. Since $N$ is generated by
$\{S^n:n=0,1,2,\cdot\cdot\cdot\}$, we conclude that
$X\equiv{M_{\theta}}$ on $N$.
\end{proof}
\begin{cor}\label{83}
Let $T\in{L(H^2)}$. For every operator
$X\in{\{S_{L(H^2)}\}^{\prime}}$ such that $X(T)\in{\{S\}^\prime}$,
there is a function $\theta^{\prime}\in{H^\infty}$ such that
$X(S^{n}T)={M_{\theta^\prime}(S^n)}$ for every
$n=0,1,2,\cdot\cdot\cdot$.
\end{cor}
\begin{proof}
 By equation \eqref{72}, $SX(T)=X(ST)$ and
 $X(S^{2}T)=SX(ST)=S^{2}X(T)$. Thus by the same way in the proof of Theorem
 \ref{71}, $X(S^{n}T)=S^{n}X(T)$ for any
 $n=0,1,2,\cdot\cdot\cdot$. Since $X(T)\in\{S\}^{\prime}$, there
 is a function $\theta^{\prime}\in{H^\infty}$ such that
 $X(T)=T_{\theta^\prime}$. By the same way in the proof of the Theorem
 \ref{71}, we can get $X(S^{n}T)={M_{\theta^\prime}(S^n)}$ for every
$n=0,1,2,\cdot\cdot\cdot$.
\end{proof}
\section{Invariant Subspaces for $S_{L(H^2)}$}\label{20}

In section \ref{21}, we have defined the shift and multiplication
operators on $L(H^2)$. In this section, we characterize closed
$S_{L(H^2)}$-invariant subspaces of $L(H^2)$.

\begin{prop}\label{25}\cite{D}
Let A and B be (bounded) operators on the Hilbert space $H$. The
following statements are equivalent:

(1) range $[A]\subset$ range $[B];$

(2) there exists a bounded operator $C$ on $H$ so that $A=BC.$
\end{prop}

\begin{thm}\label{30}
(i) If $M$ is a closed $S_{L(H^{2})}$-invariant subspace of
$L(H^{2})$, then there is an inner function
$\theta\in{H^{\infty}}$ such that
\begin{equation}\label{22}M=
T_{\theta}M^{\prime}=\{T_{\theta}T^{\prime}:T^{\prime}\in{M^{\prime}}\},\end{equation}
where $M^{\prime}$ is a closed $S_{L(H^{2})}$-invariant subspace
of $L(H^{2})$ and $H^{2}$ is generated by
$\bigcup_{T\in{M^{\prime}}}\{f\in{H^{2}}:f\in{\texttt{range of
}T}\}$.

(ii) Conversely, if $M\subset{L(H^{2})}$ satisfies equation
(\ref{22}), then $M$ is a closed invariant subspace for
$S_{L(H^{2})}$.
\end{thm}
\begin{proof}
(i) Suppose that $M_{1}$ is a closed subspace of $H^{2}$ generated
by $\bigcup_{T\in{M}}\{f\in{H^{2}}:f\in{\texttt{range of }T}\}$
and we will denote it by \[M_{1}=\bigvee\texttt{range}(M).\]

Since $M$ is a closed $S_{L(H^{2})}$-invariant subspace of
$L(H^{2})$, by the definitions of $S_{L(H^{2})}$ and $M_{1}$, we
see that $M_{1}$ is $S$-invariant. Thus, $M_{1}$ is a closed
$S$-invariant subspace of $H^{2}$.

By the Beurling's Theorem, there is an inner function
$\theta\in{H^{\infty}}$ such that \begin{equation}\label{24}
M_{1}=\theta{H^{2}}.
\end{equation}
It follows that $M_1$ is the range of $T_{\theta},$ and so for any
$T\in{M},$
\[\texttt{range }[T]\subset{M_{1}}=\texttt{range }[T_{\theta}].\]

By Proposition \ref{25}, there is an operator
$T^{\prime}\in{L(H^{2})}$ such that \[T=T_{\theta}T^{\prime}.\]

Let
\[M^{\prime}=\{T^{\prime}\in{L(H^{2})}:T=T_{\theta}T^{\prime}\texttt{ for some }T\in{M}\},\]
and $\{T_{n}\}_{n=1}^{\infty}$ be a convergent sequence in
$M^{\prime},$ such that
\begin{equation}\label{26}\lim_{n\rightarrow\infty}T_{n}=T^{\prime\prime}.\end{equation}

By definition of $M^{\prime}$, $T_{\theta}T_{n}\in{M}$ and since
\[\norm{T_{\theta}T_{n}-T_{\theta}T^{\prime\prime}}\leq\norm{\theta}_{\infty}\norm{T_{n}-T^{\prime\prime}},\]
equation (\ref{26}) implies that
$\lim_{n\rightarrow\infty}T_{\theta}T_{n}=T_{\theta}T^{\prime\prime}$.
Since $M$ is closed, $T_{\theta}T^{\prime\prime}\in{M}$, and so
$T^{\prime\prime}\in{M^{\prime}}.$ Thus $M^{\prime}$ is a closed
subspace of $L(H^{2})$.

Since $M$ is $S_{L(H^{2})}$-invariant and
$T_{\theta}T^{\prime}\in{M}$ for any $T^{\prime}\in{M^{\prime}}$,
we have
\[S_{L(H^{2})}(T_{\theta}T^{\prime})=S(T_{\theta}T^{\prime})=T_{\theta}(ST^{\prime})\in{M}.
\] It follows that \[S_{L(H^{2})}(T^{\prime})=ST^{\prime}\in{M^{\prime}}.\]
Thus $M^{\prime}$ is $S_{L(H^{2})}$-invariant.

Since $M^{\prime}$ is closed $S_{L(H^{2})}$-invariant, by the same
way as above, there is an inner function
$\theta^\prime\in{H^\infty}$ such that
\begin{equation}\label{27}\bigvee\texttt{range}(M^{\prime})=\theta^{\prime}H^{2}.\end{equation}
It follows from Proposition \ref{25} that
\[M_{1}=\theta\theta^{\prime}H^{2}.\] From equation (\ref{24}), we
have
\[\theta{H}^{2}=\theta\theta^{\prime}H^{2}.\]
By the Beurling's Theorem, $\theta^{\prime}\equiv{1}$. Thus
equation (\ref{27}) implies that
\[\bigvee\texttt{range}(M^{\prime})=H^{2}.\]

\vskip0.2cm

(ii) Suppose $M\subset{L(H^{2})}$ satisfies equation (\ref{22}).
For any $T=T_{\theta}T^{\prime}\in{M}$, i.e.
$T^{\prime}\in{M^{\prime}}$, since $M^{\prime}$ is
$S_{L(H^{2})}$-invariant, we have
\[S_{L(H^{2})}(T)=ST_{\theta}T^{\prime}=T_{\theta}ST^{\prime}=T_{\theta}(S_{L(H^{2})}(T^{\prime}))\in{T_{\theta}
M^{\prime}}=M.\] Thus $M$ is a $S_{L(H^{2})}$-invariant subspace
of $L(H^{2})$.

Let $\{T_{\theta}T_{n}\}_{n=1}^{\infty}\in{M}$, i.e.
$T_{n}\in{M^{\prime}}$, be a convergent sequence and let
\[\lim_{n\rightarrow\infty}{T}_{\theta}T_{n}=T_{0}.\]
Then for any $f\in{H^{2}}$,
\begin{equation}\label{28}{T}_{\theta}T_{n}(f)\rightarrow{T}_{0}(f)\texttt{ as
}n\rightarrow\infty\texttt{ in }H^{2}\end{equation} Since for any
$n=1,2,\cdot\cdot\cdot$ and $f\in{H^{2}}$,
$T_{\theta}{T}_{n}(f)\in{\theta{H^{2}}}$ and $\theta{H^{2}}$ is
closed in $H^{2}$, from (\ref{28}) we have
\[T_{0}f\in{\theta{H^{2}}}\texttt{ for any }f\in{H^{2}}\texttt{, i.e. range }[T_0]\subset\texttt{ range }[T_{\theta}].\]
By Proposition \ref{25}, there is $T_{1}\in{L(H^{2})}$ such that
\[T_{0}=T_{\theta}T_{1}.\]

Since
\[\norm{T_{\theta}T_{n}-T_{0}}=\norm{T_{\theta}T_{n}-T_{\theta}T_{1}}=\norm{T_{\theta}(T_{n}-T_{1})}\rightarrow
{0}\texttt{ as }n\rightarrow\infty\] and $\theta$ is inner,
\begin{equation}\label{29}\norm{T_{n}-T_{1}}\rightarrow{0}\texttt{ as }
n\rightarrow\infty.\end{equation} Since $M^{\prime}$ is closed,
$T_{1}\in{M^{\prime}}$. It follows that
\[T_{0}=T_{\theta}T_{1}\in{T_{\theta}M^{\prime}}=M.\]
Thus $M$ is a closed $S_{L(H^{2})}$-invariant subspace of
$L(H^{2})$.

\end{proof}

------------------------------------------------------------------------

\bibliographystyle{amsplain}
\bibliography{xbib}
\end{document}